\documentclass{amsart}
\usepackage{amsfonts}
\usepackage{amssymb}
\usepackage{amscd}

\newtheorem{theorem}{Theorem}[section]

\newtheorem{corollary}[theorem]{Corollary}
\newtheorem{lemma}[theorem]{Lemma}
\newtheorem{proposition}[theorem]{Proposition}

\theoremstyle{definition}
\theoremstyle{remark}
\theoremstyle{problem}
\numberwithin{equation}{section}

\def\0{\varnothing}

\def\to{\mathop{\rightarrow}\limits}

\begin{document}
\title{Group topologies coarser than the Isbell topology}
\author{Szymon Dolecki, Francis Jordan and Fr\'{e}d\'{e}ric Mynard}
\date{\today}
\maketitle

\begin{abstract}
The Isbell, compact-open and point-open topologies on the set $C(X,\mathbb{R}%
)$ of continuous real-valued maps can be represented as the dual topologies
with respect to some collections $\alpha (X)$ of compact families of open
subsets of a topological space $X$. Those $\alpha (X)$ for which addition is
jointly continuous at the zero function in $C_{\alpha }(X,{\mathbb{R}})$ are
characterized, and sufficient conditions for translations to be continuous
are found. As a result, collections $\alpha (X)$ for which $C_{\alpha }(X,{%
\mathbb{R}})$ is a topological vector space are defined canonically. The
Isbell topology coincides with this vector space topology if and only if $X$
is infraconsonant. Examples based on measure theoretic methods, that $%
C_{\alpha }(X,{\mathbb{R}})$ can be strictly finer than the compact-open
topology, are given. To our knowledge, this is the first example of a splitting group topology strictly finer than the compact-open topology. 
\end{abstract}

\section{Introduction}

The \emph{Isbell}, \emph{compact-open} and \emph{point-open} topologies on
the set $C(X,Y)$ of continuous real-valued maps from $X$ to $Y$, can be
represented as the \emph{dual topologies} with regard to some collections $%
\alpha =\alpha (X)$ of compact openly isotone families of a topological
space $X$, that is, the topology $\alpha (X,Y)$ is determined by a subbase
of open sets of the form%
\begin{equation}
\lbrack \mathcal{A},U]:=\left\{ f\in C(X,\mathbb{R}):f^{-}\left( U\right)
\in \mathcal{A}\right\} ,  \label{subbasic}
\end{equation}%
where $\varnothing \notin \mathcal{A}\in \alpha $ and $U$ are open subsets
of $Y$ (and $f^{-}(U):=\left\{ x\in X:f(x)\in U\right\} $\footnote{$f^{-}(U)$
is a shorthand for $f^{-1}(U)$}). They are dual with regard to the
collections, respectively, $\kappa (X)$ of all \emph{compact families}, $k(X)
$ of \emph{compactly generated families} and $p(X)$ \emph{finitely generated
families} on $X$. In most cases, $\alpha$ is a topology on the space $\mathcal O_X$ of open subsets of $X$. The interrelations between the topology of $X$, the topology $\alpha(X)$, and the topology $\alpha(X,Y)$ can be, to a large extent, studied at this level of generality. In particular, under mild assumptions, local properties of $\alpha(X,\mathbb R)$ at the zero function exactly correspond to local properties of $(\mathcal O_X,\alpha)$ at $X$. This greatly simplifies the study of local topological properties for the function spaces $C_\alpha(X,\mathbb R)$, provided that $C_\alpha(X,\mathbb R)$ is homogeneous.  But while $p(X,\mathbb{R})$ and $k(X,\mathbb{R})$ are
topological vector spaces for each $X$, the Isbell topology $\kappa (X,%
\mathbb{R})$ need not be even translation-invariant. In this context, it is natural to study under what conditions $\alpha(X,\mathbb R)$ is translation-invariant, or even better, a topological group. If $X$ is \emph{%
consonant} (that is, if $k(X,\$^{\ast })=\kappa (X,\$^{\ast })$, where $%
\$^{\ast }$ designs the \emph{Sierpi\'{n}ski topology}) then $k(X,\mathbb{R})
$ and $\kappa (X,\mathbb{R})$ coincide, and in particular $\kappa (X,\mathbb{%
R})$ is a group topology. In \cite{dolmyn.isbell} we characterized those
topologies $X$, for which addition is jointly continuous at the zero
function for the Isbell topology $\kappa (X,\mathbb{R})$; the class of such
topologies, called \emph{infraconsonant}, is larger than that of consonant
topologies, but we do not know if the two classes coincide in case of
completely regular topologies $X$. In this paper we prove that the Isbell
topology $\kappa (X,\mathbb{R})$ is a group topology if and only if $X$ is
infraconsonant. Therefore the problem \cite[Problem 1.1]{dolmyn.isbell} of finding a completely regular infraconsonant space that is not consonant if equivalent to the problem of finding a completely regular space $X$ such that the Isbell topology $\kappa(X,\mathbb R)$ is a group topology, but does not coincide with the compact-open topology. While we do not solve this problem, we obtain the intermediate result that there are spaces $X$ for which topologies $\alpha\subseteq\kappa(X)$ making $C_\alpha(X,\mathbb R)$ a topological group strictly finer than $C_k(X,\mathbb R)$ exist.
Indeed, for each $X$ there exists a largest
hereditary (\footnote{%
A collection $\alpha $ is \emph{hereditary }if $\mathcal{A}\downarrow A\in
\alpha $ whenever $A\in \mathcal{A}\in \alpha ,$ where $\mathcal{A}%
\downarrow A$ is defined by (\ref{eq:restrict}).}) collection $\Lambda
^{\downarrow }(X)\subseteq \kappa (X)$, for which the addition is jointly
continuous at the zero function in $\Lambda ^{\downarrow }(X,\mathbb{R})$.
It turns out that $\Lambda ^{\downarrow }(X,\mathbb{R)}$ is a vector space
topology and that a completely regular space $X$ is infraconsonant if and
only if $\Lambda ^{\downarrow }(X,\mathbb{R)=\kappa }(X,\mathbb{R)}$. Using
measure theoretic methods, we show in particular that if a completely
regular $X$ is not pre-Radon, then $k(X,\mathbb{R})$ is strictly included in 
$\Lambda ^{\downarrow }(X,\mathbb{R)}$. Therefore each completely regular space $X$ that is not pre-Radon admits a group topology on $C(X,\mathbb R)$ that is splitting (\footnote{that is, if for each space $Y$, the continuity of a map $g:Y\times X\to \mathbb R$ implies the continuity of the map $g^*:Y\to C_\alpha(X,\mathbb R)$ defined by $g^*(y)(x)=g(x,y)$; equivalently, a topology is splitting if it is coarser than the continuous convergence.}) but strictly finer than the compact-open topology. To our knowledge, no such example was known so far.

\section{Generalities}

If $\mathcal{A}$ is a family of subsets of a topological space $X$ then $%
\mathcal{O}_{X}(\mathcal{A})$ denotes the family of open subsets of $X$
containing an element of $\mathcal{A}$. In particular, if $A\subset X$ then $%
\mathcal{O}_{X}(A)$ denotes the family of open subsets of $X$ containing $A$%
. We denote by $\mathcal{O}_{X}$ the set of open subsets of $X.$

If $X$ and $Y$ are topological spaces, $C(X,Y)$ denotes the set of
continuous functions from $X$ to $Y.$ If $A\subset X$, $U\subset Y$, then $%
[A,U]:=\{f\in C(X,Y):f(A)\subset U\}$. A family $\mathcal{A}$ of subsets of $%
X$ is \emph{openly isotone} if $\mathcal{O}_{X}\left( \mathcal{A}\right) =%
\mathcal{A}$. If $\mathcal{A}$ is openly isotone and $U$ is open, then $[%
\mathcal{A},U]=\bigcup\nolimits_{A\in \mathcal{A}}[A,U]$.

If $\alpha $ is a collection of openly isotone families $\mathcal{A}$ of
open subsets of $X$, such that each open subset of $X$ belongs to an element
of $\alpha ,$ then%
\begin{equation*}
\{[\mathcal{A},U]:\mathcal{A}\in \alpha ,U\in \mathcal{O}_{Y}\}
\end{equation*}%
forms a subbase for a topology $\alpha (X,Y)$ on $C(X,Y).$ We denote the set 
$C(X,Y)$ endowed with this topology by $C_{\alpha }(X,Y)$. Note that because 
\begin{equation*}
\lbrack \mathcal{A},U]\cap \lbrack \mathcal{B},U]=[\mathcal{A}\cap \mathcal{B%
},U],
\end{equation*}%
$\alpha (X,Y)$ and $\alpha ^{\cap }(X,Y)$ coincide, where $\alpha ^{\cap }$
consists of finite intersections of the elements of $\alpha $. Therefore, we
can always assume that $\alpha $ is stable under finite intersections.

In the sequel, we will focus on the case where $\alpha $ consists of compact
families. A family $\mathcal{A=O}_{X}(\mathcal{A)}$ is \emph{compact }if
whenever $\mathcal{P}\subset \mathcal{O}_{X}$ and $\bigcup \mathcal{P}\in 
\mathcal{A}$ then there is a finite subfamily $\mathcal{P}_{0}$ of $\mathcal{%
P}$ such that $\bigcup \mathcal{P}_{0}\in \mathcal{A}$. Of course, for each
compact subset $K$ of $X,$ the family $\mathcal{O}_{X}(K)$ is compact.

We denote by $\kappa (X)$ the collection of compact families on $X$. Seen as
a family of subsets of $\mathcal{O}_{X}$ (the set of open subsets of $X$), $%
\kappa (X)$ is the set of open sets for the \emph{Scott topology}; hence
every union of compact families is compact, in particular $\bigcup_{K\in 
\mathcal{K}}\mathcal{O}_{X}(K)$ is compact if $\mathcal{K}$ is a family of
compact subsets of $X$. A topological space is called \emph{consonant }if
every compact family $\mathcal{A}$ is \emph{compactly generated}, that is,
there is a family $\mathcal{K}$ of compact sets such that $\mathcal{A}%
=\bigcup_{K\in \mathcal{K}}\mathcal{O}_{X}(K)$. Similarly, $p(X):=\{\mathcal{%
O}_{X}(F):F\in \lbrack X]^{<\omega }\}$ and $k(X):=\{{\mathcal{O}}%
(K):K\subseteq X\text{ compact }\}$ are basis for topologies on ${\mathcal{O}%
}_{X}$. Accordingly, $p(X,Y)$ is the \emph{topology of pointwise convergence}%
, $k(X,Y)$ is the \emph{compact-open topology} and $\kappa (X,Y)$ is the 
\emph{Isbell topology }on $C(X,Y).$

If $\$^{\ast }:=\left\{ \varnothing ,\left\{ 0\right\} ,\left\{ 0,1\right\}
\right\} $ the function spaces $C(X,\$^{\ast })$ can be identified with the
set of open subsets of $X$ (\footnote{%
In \cite{dolecki.pannonica}, \cite{dolmyn.isbell} and \cite{DMtransfer}, we
distinguish two homeomorphic copies $\$:=\left\{ \varnothing ,\left\{
1\right\} ,\left\{ 0,1\right\} \right\} $ and $\$^{\ast }:=\left\{
\varnothing ,\left\{ 0\right\} ,\left\{ 0,1\right\} \right\} $ of the\emph{\
Sierpi\'{n}ski topology} on $\left\{ 0,1\right\} $ and identify the function
spaces $C(X,\$)$ and $C(X,\$^{\ast })$ with the set of closed subsets of $X$
and open subsets of $X$ respectively. This is why we use $\$^{\ast }$here.}%
). In this notation, $X$ is consonant if and only if $C_{k}(X,\$^{\ast
})=C_{\kappa }(X,\$^{\ast }).$

More generally, a space $X$ is called $Z$-\emph{consonant} if $C_{\kappa
}(X,Z)=C_{k}(X,Z)$ \cite[chapter 3]{openproblems2}. \cite[Problem 62]%
{openproblems2} asks for what spaces $Z$ (other than $\$^{\ast }$) $Z$%
-consonance impies consonance. A still more general problem is, given a
collection $\alpha $ of compact families defined for each space $X,$ to
determine for what spaces $Z$,%
\begin{equation}
C_{k}(X,Z)=C_{\alpha }(X,Z)\Longleftrightarrow C_{k}(X,\$^{\ast })=C_{\alpha
}(X,\$^{\ast }).  \label{eq:alphaconsonance}
\end{equation}

The latter equality always implies the former. More generally, in view of
the definition of $\alpha (X,Z)$, if $\alpha $ and $\gamma $ are collections
(of compact families on $X$), then%
\begin{equation*}
C_{\alpha }(X,\$^{\ast })\leq C_{\gamma }(X,\$^{\ast })\Longrightarrow
C_{\alpha }(X,Z)\leq C_{\gamma }(X,Z)
\end{equation*}%
for every topological space $Z$. To show the converse implication under some
additional assumptions, recall that the \emph{restriction of }$\mathcal{A}$ 
\emph{to} $A\in \mathcal{A}$ is defined by%
\begin{equation}
\mathcal{A}\downarrow A:=\{U\in \mathcal{O}_{X}:\exists B\subseteq A\cap
U,B\in \mathcal{A}\}.  \label{eq:restrict}
\end{equation}%
\cite[Lemma 2.8]{dolmyn.isbell} shows that if $\mathcal{A}$ is a compact
family and $A\in \mathcal{A}$, then $\mathcal{A}\downarrow A$ is compact
too. A collection $\alpha $ of families of open subsets of a given set is 
\emph{hereditary }if $\mathcal{A}\downarrow A\in \alpha $ whenever $\mathcal{%
A}\in \alpha $ and $A\in \mathcal{A}$.

It was shown in \cite[Proposition 2.4]{dolmyn.isbell} that if $X$ is
completely regular and $\mathbb{R}$-consonant, then it is consonant. More
generally:

\begin{proposition}
If $\alpha ,\gamma \subseteq \kappa (X)$ are two topologies, $\alpha $ is
hereditary, $X$ is completely regular, and $C_{\alpha }(X,\mathbb{R})\leq
C_{\gamma }(X,\mathbb{R})$, then $C_{\alpha }(X,\$^{\ast })\leq C_{\gamma
}(X,\$^{\ast })$.
\end{proposition}

\begin{proof}
The neighborhood filter of an open set $A$ with respect to $\alpha
(X,\$^{\ast })$ is generated by a base of the form $\left\{ \mathcal{A}\in
\alpha :A\in \mathcal{A}\right\} $. Therefore we need show that for each $%
\mathcal{A}\in \alpha $ and each $A\in \mathcal{A},$ there exists $\mathcal{G%
}\in \gamma $ such that $\mathcal{G}\subseteq \mathcal{A}\downarrow A$. By
assumption, $\mathcal{N}_{\gamma }(\overline{0})\geq \mathcal{N}_{\alpha }(%
\overline{0})$ so that for each $\mathcal{A}\in \alpha $ and each $A\in 
\mathcal{A},$ there exists $\mathcal{G}\in \gamma $ and $r>0$ such that $%
\left[ \mathcal{G},(-r,r)\right] \subset \left[ \mathcal{A}\downarrow A,(-%
\frac{1}{2},\frac{1}{2})\right] .$ Suppose that there exists $G\in \mathcal{G%
}\setminus \left( \mathcal{A}\downarrow A\right) $, hence $X\setminus G\in
\left( \mathcal{A}\downarrow A\right) ^{\#}$. Because $X$ is completely
regular and $\mathcal{G}$ is compact there is $G_{0}\in \mathcal{G}$ and a
continuous function $f$ such that $f(G_{0})=\left\{ 0\right\} $ and $%
f(X\setminus G)=\left\{ 1\right\} ,$ by \cite[Lemma 2.5]{dolmyn.isbell}.
Then $f\in \left[ \mathcal{G},(-r,r)\right] $ but $f\notin \left[ \mathcal{A}%
\downarrow A,(-\frac{1}{2},\frac{1}{2})\right] $, because $1\in f(B)$ for
each $B\in \mathcal{A}\downarrow A$. Therefore $A\in \mathcal{G}\subseteq 
\mathcal{A}\downarrow A\subseteq \mathcal{A},$ so that $\alpha \leq \gamma .$
\end{proof}

\begin{corollary}
\label{prop:Rcoincide} If $X$ is completely regular and $\alpha \subseteq
\kappa (X)$ is hereditary, then (\ref{eq:alphaconsonance}) holds for $Z=%
\mathbb{R}$.
\end{corollary}

The \emph{grill }of\emph{\ }a family $\mathcal{A}$ of subsets of $X$ is the
family $\mathcal{A}^{\#}:=\{B\subseteq X:\forall A\in \mathcal{A}$, $A\cap
B\neq \varnothing \}$. Note that if $\mathcal{A}=\mathcal{O}(\mathcal{A)}$,
then 
\begin{equation*}
A\in \mathcal{A\Longleftrightarrow }A^{c}\notin \mathcal{A}^{\#}\text{.}
\end{equation*}

If $\mathcal{A}\in \kappa (X)$ and $C$ is a closed subset of $X$ such that $%
C\in \mathcal{A}^{\#}$ then the family%
\begin{equation*}
\mathcal{A}\vee C:=\mathcal{O}\left( \{A\cap C:A\in \mathcal{A\}}\right) ,
\end{equation*}%
called \emph{section of }$\mathcal{A}$ by $C,$ is a compact family on $X$ 
\cite{dolecki.pannonica}. A collection $\alpha $ of families of open subsets
of a given set is \emph{sectionable }if $\mathcal{A}\vee C\in \alpha $
whenever $\mathcal{A}\in \alpha $ and $C$ is a closed set in $\mathcal{A}%
^{\#}.$ It was shown in \cite[Theorem 2.9]{dolmyn.isbell} that $C_{\kappa
}(X,Z)$ is completely regular whenever $Z$ is. A simple modification of the
proof leads to the following generalization.

\begin{theorem}
\label{thm:dual_Tikh} If $Z$ is completely regular and $\alpha \subseteq
\kappa (X)$ is sectionable, then $C_{\alpha }(X,Z)$ is completely regular.
\end{theorem}

As $r[\mathcal{A},U]=[\mathcal{A},rU]$ for all $r\neq 0,$ it is immediate
that inversion for $+$ is always continuous in $C_{\alpha }(X,\mathbb{R}).$
More generally, the proof of the joint continuity of scalar multiplication
in $C_{\kappa }(X,\mathbb{R)}$ \cite[Proposition 2.10]{dolmyn.isbell} can be
adapted to the effect that:

\begin{proposition}
\label{prop:scalarmultIsbell} If $\alpha \subseteq \kappa (X)$ is
hereditary, then multiplication by scalars is jointly continuous for $%
C_{\alpha }(X,\mathbb{R})$.
\end{proposition}

\begin{corollary}
Let $\alpha \subseteq \kappa (X)$ be hereditary. If $C_{\alpha }(X,\mathbb{R}%
)$ is a topological group then it is a topological vector space.
\end{corollary}

\section{Self-joinable collections and joint continuity of addition at the
zero function}

As usual, if $A$ and $B$ are subsets of an additive group, $A+B:=\{a+b:a\in
A,b\in B\}$ and if $\mathcal{A}$ and $\mathcal{B}$ are two families of
subsets, $\mathcal{A}+\mathcal{B}:=\{A+B:A\in \mathcal{A},B\in \mathcal{B}\}$%
.

As we have mentioned, a topology on an additive group is a group topology if
and only if inversion and translations are continuous, and $\mathcal{N}(o)+%
\mathcal{N}(o)\geq \mathcal{N}(o)$, where $o$ is the neutral element. First,
we investigate the latter property, that is, 
\begin{equation}
\mathcal{N}_{\alpha }(\overline{0})+\mathcal{N}_{\alpha }(\overline{0})\geq 
\mathcal{N}_{\alpha }(\overline{0}),
\end{equation}%
for the space $C_{\alpha }(X,\mathbb{R}),$ where $\overline{0}$ denotes the
zero function.

If $\alpha $ and $\gamma $ are two subsets of $\kappa (X),$ we say that $%
\alpha $ is $\gamma $-\emph{joinable }if for every $\mathcal{A}\in \alpha ,$
there is $\mathcal{G}\in \gamma $ such that ${\mathcal{G}}\vee {\mathcal{G}}%
\subseteq {\mathcal{A}}$, where%
\begin{equation*}
{\mathcal{G}}\vee {\mathcal{G}}:=\{G_{1}\cap G_{2}:G_{1},G_{2}\in \mathcal{G}%
\}.
\end{equation*}

A subset $\alpha $ of $\kappa (X)$ is \emph{self-joinable} if it is $\alpha $%
-joinable. A family $\mathcal{A}$ is called \emph{joinable }if $\{\mathcal{A}%
\}$ is $\kappa (X)$-joinable. In \cite{dolmyn.isbell}, a space $X$ is called 
\emph{infraconsonant} if every compact family is joinable, that is, $\kappa
(X)$ is self-joinable. \cite[Theorem 3.1]{dolmyn.isbell} shows that among
completely regular spaces $X,$ 
\begin{equation*}
\mathcal{N}_{\kappa }(\overline{0})+\mathcal{N}_{\kappa }(\overline{0})\geq 
\mathcal{N}_{\kappa }(\overline{0})
\end{equation*}%
if and only if $X$ is infraconsonant. More generally,

\begin{theorem}
\label{th:at0} Let $X$ be a completely regular space. Then $\alpha \subseteq
\kappa (X)$ is self-joinable if and only if 
\begin{equation}
{\mathcal{N}}_{\alpha }(\overline{0})+{\mathcal{N}}_{\alpha }(\overline{0}%
)\geq {\mathcal{N}}_{\alpha }(\overline{0}).  \label{eq:at0b}
\end{equation}
\end{theorem}

\begin{proof}
Let $\mathcal{A\in }\alpha $ and $V\in \mathcal{N}_{\mathbb{R}}(0).$ Because 
$\alpha $ is self-joinable, there exist$\mathcal{\ }$a compact family $%
\mathcal{B}$ in $\alpha $ such that $\mathcal{B}\vee \mathcal{B}\subseteq 
\mathcal{A}$. If $W\in \mathcal{N}_{\mathbb{R}}(0)$ such that $W+W\subseteq
V $, then $[\mathcal{B},W]+[\mathcal{B},W]\subseteq \lbrack \mathcal{A},V]$,
which proves (\ref{eq:at0b})$.$

Conversely, assume that $\alpha =\alpha ^{\cap }$ is not self-joinable. Let $%
\mathcal{A}$ be a family of $\alpha $ such that $\mathcal{B\vee B\nsubseteq A%
}$ for every $\mathcal{B}\in \alpha $. Note that $\mathcal{B\vee C\nsubseteq
A}$ for every pair of families $\mathcal{B}$ and $\mathcal{C}$ in $\alpha $
for otherwise $\mathcal{D}=\mathcal{B\cap C}$ would be a family of $\alpha $
such that $\mathcal{D\vee D\subseteq A}$. Let $V=\left( -\frac{1}{2},\frac{1%
}{2}\right) .$ We claim that for any pair $(\mathcal{B},\mathcal{C)}\in
\alpha ^{2}$ and any pair $(U,W)$ of $\mathbb{R}$-neighborhood of $0,$ $[%
\mathcal{B},U]+[\mathcal{C},W]\nsubseteq \lbrack \mathcal{A},V]$. Indeed,
there exist $B\in \mathcal{B}$ and $C\in \mathcal{C}$ such that $B\cap
C\notin \mathcal{A}$. Then $B^{c}\cup C^{c}\in \mathcal{A}^{\#}$. Moreover, $%
B^{c}\notin \mathcal{B}^{\#}$ so that by \cite[Lemma 2.5]{dolmyn.isbell},
there exist $B_{1}\in \mathcal{B}$ and $f\in C(X,\mathbb{R})$ such that $%
f(B_{1})=\{0\}$ and $f(B^{c})=\{1\}.$ Similarly, $C^{c}\notin \mathcal{C}$
so that there exist $C_{1}\in \mathcal{C}$ and $g\in C(X,\mathbb{R})$ such
that $g(C_{1})=\{0\}$ and $g(C^{c})=\{1\}.$ Then $f+g\in \lbrack \mathcal{B}%
,U]+[\mathcal{C},W]$ but $1\in (f+g)(A)$ for all $A\in \mathcal{A}$ so that $%
f+g\notin \lbrack \mathcal{A},V]$.
\end{proof}

Note that the collection $k(X)$ (of compactly generated families) is
self-joinable, and that a union of self-joinable collections is
self-joinable. Therefore, there is a largest self-joinable subset $\Lambda
(X)$ of $\kappa (X)$. If $\alpha $ is self-joinable, so is $\alpha ^{\cap }$%
. Therefore $\Lambda (X)$ is stable for finite intersections. In fact, $%
\Lambda (X)$ is a topology on $C(X,\$^{\ast })$ and%
\begin{equation}
k(X)\subseteq \Lambda (X)\subseteq \kappa (X).  \label{join-ineq}
\end{equation}

\begin{corollary}
\label{cor:join} Let $X$ be completely regular. The largest subcollection $%
\alpha $ of $\kappa (X)$, for which (\ref{eq:at0b}) holds is $\alpha
=\Lambda (X)$. In particular, a completely regular space $X$ is
infraconsonant if and only if $\kappa (X)=\Lambda (X).$
\end{corollary}

The collection $\Lambda (X)$ is sectionable (\footnote{%
Indeed, if $\alpha $ is self-joinable, so is $\left\{ \mathcal{A}\vee C:%
\mathcal{A}\in \alpha ,C\in \mathcal{A}^{\#},C\text{ closed}\right\} $
because if $\mathcal{A}\in \alpha $, there is $\mathcal{B}\in \alpha $ such
that $\mathcal{B}\vee \mathcal{B}\subset \mathcal{A}\subset \mathcal{A}\vee
C $. By maximality, $\Lambda (X)$ is sectionable.}), but in general it is
not hereditary. We will construct now a largest hereditary collection of
compact families for which (\ref{eq:at0b}) holds.

If $\alpha $ and $\gamma $ are two subsets of $\kappa (X),$ we say that $%
\alpha $ is \emph{hereditarily }$\gamma $-\emph{joinable }if for every $%
\mathcal{A}\in \alpha ,$ and every $A\in \mathcal{A}$, there is $\mathcal{G}%
\in \gamma $ such that $A\in \mathcal{G}$ and ${\mathcal{G}}\vee {\mathcal{G}%
}\subseteq {\mathcal{A}}$. A subset $\alpha $ of $\kappa (X)$ is \emph{%
hereditarily} \emph{self-joinable} if it is hereditarily $\alpha $-joinable.
A family $\mathcal{A}$ is called \emph{hereditarily} \emph{joinable }if $\{%
\mathcal{A}\}$ is hereditarily $\kappa (X)$-joinable. There exists a largest
hereditarily self-joinable subset $\Lambda ^{\downarrow }(X)$ of $\kappa
(X). $ Notice that $\Lambda ^{\downarrow }(X)$ is also the largest
self-joinable and hereditary collection of compact families, and that $%
\Lambda ^{\downarrow }(X)$ is sectionable.

\begin{corollary}
Let $X$ be completely regular. The largest hereditary subcollection $\alpha $
of $\kappa (X)$, for which (\ref{eq:at0b}) holds is $\alpha =\Lambda
^{\downarrow }(X)$. In particular, a completely regular space $X$ is
infraconsonant if and only if $\kappa (X)=\Lambda ^{\downarrow }(X).$
\end{corollary}

Of course, $\Lambda ^{\downarrow }(X)$ is a topology, and $\Lambda
^{\downarrow }(X)\subseteq \Lambda (X)$. The inclusion can be strict. In
fact, we have:

\begin{proposition}
\label{prop:heredinfra} A regular space $X$ is infraconsonant if and only if 
$\Lambda ^{\downarrow }(X)=\Lambda (X)$ if and only if $\kappa (X)=\Lambda
^{\downarrow }(X).$
\end{proposition}

\begin{proof}
If $X$ is not infraconsonant, there is a non-joinable family $\mathcal{A}$
on $X$. For any $x\in X\setminus \bigcap \mathcal{A},$ the family $\mathcal{O%
}(x)\cup \mathcal{A}$ belongs to $\Lambda (X)$ but not to $\Lambda
^{\downarrow }(X)$. If\ $X$ is regular and infraconsonant, then by \cite[%
Lemma 3.2]{dolmyn.isbell} $\kappa (X)=\Lambda ^{\downarrow }(X)$. Finally if 
$\kappa (X)=\Lambda ^{\downarrow }(X)$ then $\Lambda ^{\downarrow
}(X)=\Lambda (X)$, because $\Lambda (X)$ is between $\Lambda ^{\downarrow
}(X)$ and $\kappa (X)$.
\end{proof}

Examples of non-infraconsonant spaces are provided in \cite{dolmyn.isbell},
so that both inclusions in (\ref{join-ineq}) can be strict simultaneously.

\section{Self-splittable collections and continuity of translations}

A collection $\alpha $ is $\gamma $\emph{-splittable} if for every $\mathcal{%
A}\in \alpha $ and for every open subsets $U_{1}$ and $U_{2}$ of $X$ such
that $U_{1}\cup U_{2}\in {\mathcal{A}}$ there exist families ${\mathcal{G}}%
_{i}={\mathcal{G}}_{i}\downarrow U_{i}$ in $\gamma ,$ $i\in \{1,2\},$ such
that ${\mathcal{G}}_{1}\cap \mathcal{G}_{2}\subseteq {\mathcal{A}}$. A
collection $\alpha $ is \emph{self-splittable }if it is $\alpha $%
-splittable. A compact family $\mathcal{A}$ on $X$ is \emph{splittable }if $%
\left\{ \mathcal{A}\right\} $ is $\kappa (X)$-splittable. An immediate
induction shows that if $\alpha $ is self-splittable, then for every $%
\mathcal{A}\in \alpha $ and every finite collection $\{U_{1},\ldots ,U_{n}\}$
of open subsets of $X$ such that $\bigcup_{i=1}^{n}U_{i}\in \mathcal{A},$
there are families $\mathcal{C}_{i}\in \alpha $ with $\mathcal{C}_{i}=%
\mathcal{C}_{i}\downarrow U_{i}$ such that $\bigcap_{i=1}^{n}\mathcal{C}%
_{i}\subseteq \mathcal{A}$.

In \cite{francis.split}, F. Jordan calls a topological space \emph{compactly
splittable }if every compact family is splittable. A topological space with
at most one non-isolated point is said to be \emph{prime.} A modification of
the proof of \cite[Theorem 18]{francis.coincide} shows:

\begin{proposition}
\label{prop:prime} Prime spaces are compactly splittable.
\end{proposition}

It follows from \cite[Theorem 2]{francis.split} and \cite[Corollary 4.2]%
{dolmyn.isbell} that translations are continuous in the Isbell topological
space $C_{\kappa }(X,\mathbb{R})$ if $X$ is compactly splittable. More
generally, we have:

\begin{proposition}
\label{prop:translations} If $\alpha \subseteq \kappa (X)$ is
self-splittable, then translations are continuous for $C_{\alpha }(X,{%
\mathbb{R}})$.
\end{proposition}

\begin{proof}
We show continuity of the translation by $f_{0}$ at $g_{0}.$ Let ${\mathcal{A%
}}\in \alpha $ and $U\in {\mathcal{O}}_{{\mathbb{R}}}$ such that $%
f_{0}+g_{0}\in \lbrack {\mathcal{A}},U]$. There is $A_{0}\in {\mathcal{A}}$
such that $(f_{0}+g_{0})(A_{0})\subseteq U$. For each $x\in A_{0}$, there
exists $V_{x}=-V_{x}\in {\mathcal{O}}_{{\mathbb{R}}}(0)$ such that $%
f_{0}(x)+g_{0}(x)+2V_{x}\subseteq U$. Moreover, by continuity of $f_{0}$ and 
$g_{0}$, there is $W_{x}\in {\mathcal{O}}_{X}(x)$ such that $%
f_{0}(W_{x})\subseteq f_{0}(x)+V_{x}$ and $g_{0}(W_{x})\subseteq
g_{0}(x)+V_{x}$. As $\bigcup_{x\in A_{0}}W_{x}\in {\mathcal{A}}$ and ${%
\mathcal{A}}$ is compact, there is a finite subset $F$ of $A_{0}$ such that $%
W:=\bigcup_{x\in F}W_{x}\in {\mathcal{A}}$. Because $\alpha $ is
self-splittable, there exists, for each $x\in F$, a compact family ${%
\mathcal{C}}_{x}={\mathcal{C}}_{x}\downarrow W_{x}$ of $\alpha $ such that $%
\bigcap_{x\in F}{\mathcal{C}}_{x}\subseteq {\mathcal{A}}$. Note that by
construction $g_{0}\in \bigcap_{x\in F}[{\mathcal{C}}_{x},g_{0}(x)+V_{x}]$.
Moreover, if $g\in \bigcap_{x\in F}[{\mathcal{C}}_{x},g_{0}(x)+V_{x}]$, then
for each $x\in F$ there is $C_{x}\in {\mathcal{C}}_{x}$, $C_{x}\subseteq
W_{x}$, such that $g(C_{x})\subseteq g_{0}(x)+V_{x}$. If $y\in \bigcup_{x\in
F}C_{x}\in {\mathcal{A}}$, then $y\in C_{x}$ for some $x$ and 
\begin{equation*}
f_{0}(y)+g(y)\in f_{0}(x)+V_{x}+g_{0}(x)+V_{x}\subseteq U,
\end{equation*}%
so that $f_{0}+\bigcap_{x\in F}[{\mathcal{C}}_{x},g_{0}(x)+V_{x}]\subseteq
\lbrack {\mathcal{A}},U]$.
\end{proof}

\cite[Example 4.9]{dolecki.pannonica} shows that even on a compactly
splittable space (like the Arens space), there exists $\alpha \subseteq
\kappa (X)$ such that translations are not continuous for $C_{\alpha }(X,%
\mathbb{R}),$ so that $\alpha $ is not self-splittable.

Note that in a regular space $X$, the collection $k(X)$ is self-splittable,
and a union of self-splittable collections is self-splittable. Therefore,
there is a largest self-splittable subset $\Sigma (X)$ of $\kappa (X)$. If $%
\alpha $ is self-splittable, so is $\alpha ^{\cap }$, hence that $\Sigma (X)$
is a topology on $C(X,\$^{\ast })$, and%
\begin{equation}
k(X)\subseteq \Sigma (X)\subseteq \kappa (X).  \label{split-ineq}
\end{equation}

Moreover, $\Sigma (X)$ is clearly hereditary, and sectionable (\footnote{%
Indeed, if $\alpha $ is self-splittable, so is $\left\{ \mathcal{A}\vee C:%
\mathcal{A}\in \alpha ,C\in \mathcal{A}^{\#},C\text{ closed}\right\} $.
Indeed, if $\bigcup_{i=1}^{n}U_{i}\in \mathcal{A}\vee C,$ then $%
\bigcup_{i=1}^{n}U_{i}\cup C^{c}\in \mathcal{A},$ so that there are families 
$\mathcal{C}_{i}=\mathcal{C}_{i}\downarrow U_{i}$ and $\mathcal{C}_{c}=%
\mathcal{C}_{c}\downarrow C^{c}$ in $\alpha $ such that $\bigcap \mathcal{C}%
_{i}\cap \mathcal{C}_{c}\subseteq \mathcal{A}$. Therefore%
\begin{equation*}
\bigcap_{i}\mathcal{C}_{i}\vee C=\left( \bigcap_{i}\mathcal{C}_{i}\cap 
\mathcal{C}_{c}\right) \vee C\subseteq \mathcal{A}\vee C.
\end{equation*}%
By maximality, $\Lambda (X)$ is sectionable.}).

Both inequalities in (\ref{split-ineq}) can be strict. Examples of
non-compactly splittable spaces are provided in \cite{francis.coincide}, so
that $\Sigma (X)$ can be strictly included in $\kappa (X)$. On the other
hand, in view of Proposition \ref{prop:prime}, if $X$ is prime and not
consonant (e.g., the Arens space), then $k(X)$ is strictly included in $%
\Sigma (X)$.

\begin{theorem}
\label{thm:1}Let $X$ be regular. If $\alpha \subseteq \kappa (X)$ is
self-joinable, hereditary, and sectionable, then $\alpha $ is
self-splittable.
\end{theorem}

\begin{proof}
By way of contradiction, assume that $\alpha $ is not self-splittable. Let $%
\alpha _{1}=\alpha \cup \{\mathcal{O}(X)\}$. Using that $\emptyset \in 
\mathcal{O}(X)$, one can easily show that $\alpha _{1}$ is self-joinable,
hereditary, and sectionable. It is also easy to check that $\alpha _{1}$ is
not self-splittable.

Let $\mathcal{C}\in \alpha _{1}$ witness that $\alpha _{1}$ is not
self-splittable. There exist $\mathcal{C}\in \alpha _{1}$ and open sets $%
U_{1},U_{2}$ such that $U_{1}\cup U_{2}\in \mathcal{C}$, but for any $%
\mathcal{B}_{1},\mathcal{B}_{2}\in \alpha _{1}$ with $U_{1}\in \mathcal{B}%
_{1}$ and $U_{2}\in \mathcal{B}_{2}$ we have $\mathcal{B}_{1}\cap \mathcal{B}%
_{2}{\nsubseteq }\mathcal{C}$.

By regularity, there is an open cover $\mathcal{V}$ of $U_{1}\cup U_{2}$
such that $\mathrm{cl}(V)\subseteq U_{1}$ or $\mathrm{cl}(V)\subseteq U_{2}$ for
every $V\in \mathcal{V}$. Since $U_{1}\cup U_{2}\in \mathcal{C}$, there
exist a finite $\mathcal{V}_{1}\subseteq \mathcal{V}$ such that $\bigcup 
\mathcal{V}_{1}\in \mathcal{C}$. Let $W_{1}=\bigcup \{V\in \mathcal{V}%
_{1}\colon \mathrm{cl}(V)\subseteq U_{1}\}$ and $W_{2}=\bigcup \{V\in \mathcal{%
V}_{1}\colon \mathrm{cl}(V)\subseteq U_{2}\}$. Notice that $W_{1}\cup W_{2}\in 
\mathcal{C}$, $\mathrm{cl}(W_{1})\subseteq U_{1}$, and $\mathrm{cl}%
(W_{2})\subseteq U_{2}$. Let $\mathcal{C}_{1}=\mathcal{C}\downarrow
(W_{1}\cup W_{2})$. Since $\alpha _{1}$ is hereditary, $\mathcal{C}_{1}\in
\alpha _{1}$. By self-joinability of $\alpha _{1}$ there is a $\mathcal{D}%
\in \alpha _{1}$ such that $\mathcal{D}\bigvee \mathcal{D}\subseteq \mathcal{%
C}_{1}$. Notice that $\mathcal{D}\subseteq \mathcal{C}_{1}$.

Suppose there is a $D\in \mathcal{D}$ such that $D\cap W_{1}=\emptyset $.
Since $D\in \mathcal{D}\subseteq \mathcal{C}_{1}$, there is an $E\in 
\mathcal{C}$ such that $E\subseteq D\cap (W_{1}\cup W_{2})$. Notice that $%
E\subseteq W_{2}\subseteq U_{2}$. So, $U_{2}\in \mathcal{C}$. Since $\alpha
_{1}$ is hereditary, $\mathcal{C}\downarrow U_{2}\in \alpha _{1}$. Notice
that $U_{2}\in \mathcal{C}\downarrow U_{2}\in \alpha _{1}$, $U_{1}\in 
\mathcal{O}(X)\in \alpha _{1}$, and $\mathcal{C}\downarrow U_{2}\cap 
\mathcal{O}(X)=\mathcal{C}\downarrow U_{2}\subseteq \mathcal{C}$, which
contradicts our choice of $U_{1}$ and $U_{2}$. So, we may assume that $%
\mathcal{D}\#W_{1}$. Similarly, we may assume that $\mathcal{D}\#W_{2}$.

For each $i\in \{1,2\}$ let $\mathcal{D}_{i}=\mathcal{D}\bigvee \mathrm{cl}%
(W_{i})$. Since $\alpha _{1}$ is sectionable, $\mathcal{D}_{1},\mathcal{D}%
_{2}\in \alpha _{1}$. Notice that $U_{i}\in \mathcal{D}_{i}$ for every $i$.
Let $P\in \mathcal{D}_{1}\cap \mathcal{D}_{2}$. For every $i\in \{1,2\}$
there exist $D_{i}\in \mathcal{D}$ such that $D_{i}\cap \mathrm{cl}%
(W_{i})\subseteq P$. Since $D_{1},D_{2}\in \mathcal{D}$, $D_{1}\cap D_{2}\in 
\mathcal{C}_{1}$. There is an $E\in \mathcal{C}$ such that $E\subseteq
D_{1}\cap D_{2}\cap (W_{1}\cup W_{2})$. Now, 
\begin{equation*}
E\subseteq (D_{1}\cap D_{2})\cap \mathrm{cl}(W_{1}\cup W_{2})\subseteq
(D_{1}\cap \mathrm{cl}(W_{1}))\cup (D_{2}\cap \mathrm{cl}(W_{2}))\subseteq P.
\end{equation*}%
So, $P\in \mathcal{C}$. Thus, $\mathcal{D}_{1}\cap \mathcal{D}_{2}\subseteq 
\mathcal{C}$, contradicting our choice of $U_{1}$ and $U_{2}$.
\end{proof}

\begin{corollary}
\label{cor:1} Let $X$ be regular. If $X$ is infraconsonant then $X$ is
compactly splittable.
\end{corollary}

\begin{proof}
If $X$ is infraconsonant, then $\kappa (X)$ is self-joinable. Since $\kappa
(X)$ is hereditary and sectionable, $\kappa (X)$ is self-splittable, by
Theorem~\ref{thm:1}. Since $\kappa (X)$ is self-splittable, $X$ is compactly
splittable.
\end{proof}

Note that Corollary \ref{cor:1} provides a negative answer to \cite[Problem
1.2]{dolmyn.isbell}. The converse of Corollary \ref{cor:1} is not true. For
instance, the Arens space is compactly splittable because it is prime, but
it is not infraconsonant \cite[Theorem 3.6]{dolmyn.isbell}.

In the diagram below, $X$ is a regular space, and arrows represent
inclusions. We have already justified that all of these inclusions may be
strict, except for $k(X)\subseteq \Lambda ^{\downarrow }(X)$. We will see in
the next section that it may be strict.

\begin{center}
\begin{picture}(210,100)

\put(60,50){\makebox(0,0){$k(X)$}}
\put(110,50){\makebox(0,0){$\Lambda^{\downarrow}(X)$}}
\put(150,90){\makebox(0,0){$\Sigma(X)$} }
\put(150,10){\makebox(0,0){$\Lambda(X)$}}
\put(190,50){\makebox(0,0){$\kappa(X)$} }

\put(80,50){\vector(1,0){15}}
\put(115,65){\vector(1,1){20}}
\put(115,35){\vector(1,-1){20}}
\put(165,80){\vector(1,-1){20}}
\put(165,20){\vector(1,1){20}}
\end{picture}
\end{center}

In view of Proposition \ref{prop:scalarmultIsbell}, Theorem \ref{th:at0},
Proposition \ref{prop:translations} and Theorem \ref{thm:1}, we obtain:

\begin{corollary}
\label{cor:group} Let $X$ be completely regular. Then $C_{\Lambda
^{\downarrow }}(X,\mathbb{R)}$ is a topological vector space.
\end{corollary}

As a consequence, we can extend \cite[Theorem 5.3]{dolmyn.isbell} from prime
spaces to general completely regular spaces to the effect that:

\begin{corollary}
\label{cor:infragroup} Let $X$ be completely regular. The following are
equivalent:

\begin{enumerate}
\item $X$ is infraconsonant;

\item $\kappa (X)=\Lambda (X);$

\item $\kappa (X)=\Lambda ^{\downarrow }(X);$

\item $\Lambda (X)=\Lambda ^{\downarrow }(X);$

\item $C_{\kappa }(X,\mathbb{R})$ is a topological vector space;

\item $C_{\kappa }(X,\mathbb{R})$ is a topological group;

\item $\mathcal{N}_{\kappa }(\overline{0})+\mathcal{N}_{\kappa }(\overline{0}%
)\geq \mathcal{N}_{\kappa }(\overline{0});$

\item $\cap :C_{\kappa }(X,\$^{\ast })\times C_{\kappa }(X,\$^{\ast
})\rightarrow C_{\kappa }(X,\$^{\ast })$ is jointly continuous.
\end{enumerate}
\end{corollary}

\begin{proof}
Equivalences between (1) through (7) follow immediately from Theorem \ref%
{th:at0}, Corollary \ref{cor:join} and Proposition \ref{prop:heredinfra}.
The equivalence with (8) follows from \cite[Proposition 3.3]{dolmyn.isbell}.
\end{proof}

\section{A vector space topology strictly finer than the compact-open
topology}

A finite measure $\mu $ is called $\tau $-\emph{additive }if for every
family $\mathcal{P}\subset \mathcal{O}_{X},$ and for every $\varepsilon >0$
there is a finite subfamily $\mathcal{P}_{\varepsilon }\subset \mathcal{P}$
such that $\mu (\bigcup \mathcal{P}_{\varepsilon })\geq \mu (\bigcup 
\mathcal{P})-\varepsilon .$ Hence, if $\mu $ is a $\tau $-additive measure
on $X,$ then for each $r>0,$ the family 
\begin{equation*}
\mathcal{M}_{r}^{\mu }:=\left\{ O\in \mathcal{O}_{X}:\mu (O)>r\right\}
\end{equation*}%
is compact. A topological space $X$ is called \emph{pre-Radon }if every
finite $\tau $-additive measure $\mu $ on $X$ is a Radon measure, that is, $%
\mu (B)=\sup \left\{ \mu (K):K\subset B,K\text{ compact}\right\} $ for each
Borel subset $B$ of $X$.

\begin{lemma}
\label{lem:familymeasure} Let $\mu $ be a $\tau $-additive finite measure on
a space $X$. Then $\gamma _{\mu }:=\{\mathcal{M}_{r}^{\mu }\downarrow A:A\in 
\mathcal{M}_{r}^{\mu },r>0\}$ is self-splittable and hereditarily
self-joinable.
\end{lemma}

\begin{proof}
\emph{Proof of 1. }$\gamma _{\mu }$ is hereditary and self-joinable (hence
hereditarily self-joinable) because if $U\in \mathcal{M}_{r}^{\mu }$ and $m=%
\frac{r+\mu (U)}{2}$ then%
\begin{equation}
\left( \mathcal{M}_{m}^{\mu }\downarrow U\right) \vee \left( \mathcal{M}%
_{m}^{\mu }\downarrow U\right) \subseteq \mathcal{M}_{r}^{\mu }\downarrow U.
\label{eq:selfjoin}
\end{equation}

Indeed, if $O_{1}$ and $O_{2}$ are elements of $\mathcal{M}_{m}^{\mu
}\downarrow U,$ we can assume that $\mu (O_{1}\cup O_{2})\leq \mu (U)$ so
that 
\begin{equation*}
\mu (O_{1}\cap O_{2})=\mu (O_{1})+\mu (O_{2})-\mu (O_{1}\cup O_{2})>2m-\mu
(U)=r.
\end{equation*}

Self-splittability follows from the fact (\footnote{%
Note that $\gamma _{\mu }$ is not sectionable, so that Theorem \ref{thm:1}
is not sufficient to deduce self-splittability.}) that if $U_{1}\cup
U_{2}\in {\mathcal{M}}_{r}$, for $d:=\min (\mu (U_{1}),\mu (U_{2}),\mu
(U_{1}\cup U_{2})-r)>0$, $m_{1}:=\mu (U_{1})-\frac{d}{2}$, $m_{2}:=\mu
(U_{2})-\frac{d}{2}$, we have 
\begin{equation}
\left( {\mathcal{M}}_{m_{1}}^{\mu }\downarrow U_{1}\right) \cap \left( {%
\mathcal{M}}_{m_{2}}^{\mu }\downarrow U_{2}\right) \subseteq {\mathcal{M}}%
_{r}.  \label{eq:selfsplit}
\end{equation}%
Indeed, if $A_{i}\in {\mathcal{M}}_{m_{i}}^{\mu }\downarrow U_{i}$ for $i\in
\{1,2\}$, then 
\begin{eqnarray*}
&&\mu (A_{1}\cup A_{2})=\mu (A_{1})+\mu (A_{2})-\mu (A_{1}\cap A_{2}) \\
&\geq &\mu (A_{1})+\mu (A_{2})-\mu (U_{1}\cap U_{2}) \\
&>&\mu (U_{1})+\mu (U_{2})-\mu (U_{1}\cap U_{2})-d \\
&>&r.
\end{eqnarray*}
\end{proof}

\begin{theorem}
\label{th:preradon} If $X$ is a (Hausdorff) completely regular space but is
not pre-Radon, then 
\begin{equation*}
C_{\Lambda ^{\downarrow }}(X,{\mathbb{R}})>C_{k}(X,{\mathbb{R}}).
\end{equation*}
\end{theorem}

\begin{proof}
The proof of \cite[Proposition 3.1]{bouziad.borel} shows that if $X$ is a
Hausdorff non pre-Radon space, then there is a $\tau $-additive finite
measure $\mu $ and an $r>0$ such that ${\mathcal{M}}_{r}$ is compact but not
compactly generated. In view of Lemma \ref{lem:familymeasure}, $\gamma _{\mu
}\subseteq \Lambda ^{\downarrow }(X)$ so that $C_{k}(X,\$^{\ast
})<C_{\Lambda ^{\downarrow }}(X,\$^{\ast })$. In view of Proposition \ref%
{prop:Rcoincide}, $C_{k}(X,\mathbb{R})<C_{\Lambda ^{\downarrow }}(X,\mathbb{R%
})$ because $\Lambda ^{\downarrow }(X)$ is hereditary.
\end{proof}

Note that we have shown that the inclusion $k(X)\subseteq \Lambda
^{\downarrow }(X)$ may be strict.

For instance, if $X$ is the Sorgenfrey line, which is not pre-Radon, then $%
C_{\Lambda ^{\downarrow }}(X,{\mathbb{R}})$ is a topological vector space
and is strictly finer than the compact-open topology.


\begin{thebibliography}{99}
\bibitem{alleche.calbrix} B. Alleche and J. Calbrix, \emph{On the
coincidence of the upper Kuratowski topology with the cocompact topology},
Topology Appl., \textbf{93}(1999), 207-218.

\bibitem{bouziad.borel} A. Bouziad, \emph{Borel measures in consonant spaces}%
, Topology Appl., \textbf{70} (1996), 125-132.

\bibitem{dolecki.pannonica} S.~Dolecki, \emph{Properties transfer between
topologies on function spaces, hyperspaces and underlying spaces},
Mathematica Pannonica, \textbf{19}(2) (2008), 243-262.

\bibitem{DGL.kur} S.~Dolecki, G.~H. Greco, and A.~Lechicki, \emph{When do
the upper {K}uratowski topology (homeomorphically, {S}cott topology) and the
cocompact topology coincide?}, Trans. Amer. Math. Soc. \textbf{\textbf{347}}
(1995), 2869--2884.

\bibitem{dolmyn.isbell} S. Dolecki and F. Mynard, \emph{When is the Isbell
topology a group topology?}, to appear in Topology Appl.

\bibitem{DMtransfer} S. Dolecki and F. Mynard. \emph{Relations between
function spaces, hyperspaces and underlying spaces.} in preparation.

\bibitem{Isbell_a} J.R. Isbell, \emph{Function spaces and adjoints}, Math.
Scandinavica \textbf{\textbf{36}} (1975), 317--339.

\bibitem{Isbell_b} J.R. Isbell, \emph{Meet-continuous lattices}, Symposia
Mathematica \textbf{\textbf{16}} (1975), 41--54.

\bibitem{francis.coincide} F. Jordan, \emph{Coincidence of function space
topologies}, submitted, 2008.

\bibitem{francis.split} F. Jordan, \emph{More on coincidence of function
space topologies}, in preparation.

\bibitem{openproblems2} Elliot Pearl (ed.), \emph{Open problems in topology {%
II}}, Elsevier, 2007.
\end{thebibliography}
\end{document}